\documentclass[11pt]{amsart}
\usepackage{graphicx}
\newtheorem{thm}{Theorem}[section]

\newtheorem{lem}[thm]{Lemma}

\theoremstyle{definition}

\theoremstyle{remark}

\numberwithin{equation}{section}

\begin{document}

\title[Optimal control of a large dam]{Optimal control of a large dam}%
\author{Vyacheslav M. Abramov}%
\address{School of Mathematical Sciences, Monash University, Building 28M,
Clayton Campus, Clayton, Victoria 3800}%
\email{vyacheslav.abramov@sci.monash.edu.au}%


\begin{abstract}
A large dam model is an object of study of this paper. The
parameters $L^{lower}$ and $L^{upper}$ are its lower and upper
levels, $L=L^{upper}-L^{lower}$ is large, and if a current level
of water is between these bounds, then the dam is assumed to be in
normal state. Passage one or other bound leads to damage. Let
$J_1$ $(J_2)$ denote the damage cost of crossing the lower (upper)
level. It is assumed that input stream of water is described by a
Poisson process, while the output stream is state-dependent (the
exact formulation of the problem is given in the paper). Let $L_t$
denote the dam level at time $t$, and let
$p_1=\lim_{t\to\infty}\mathbf{P}\{L_t= L^{lower}\}$,
$p_2=\lim_{t\to\infty}\mathbf{P}\{L_t> L^{upper}\}$ exist. The
long-run average cost $J=p_1J_1+p_2J_2$ is a performance measure.
The aim of the paper is to choose
 the parameter of output stream (exactly specified in the paper)
minimizing $J$.

\end{abstract}
\maketitle
\section{\bf Introduction}
A large dam model is an object of study of this paper. The
parameters $L^{lower}$ and $L^{upper}$ are lower and upper levels
of the dam, and if a current level of water is between these
bounds, then the dam is assumed to be in normal state. The reason
that a dam is called \emph{large} is that the difference
$L=L^{upper}-L^{lower}$ is large. This assumption enables us to
use asymptotic analysis as $L\to\infty$ and obtain much simpler
representations for the desired characteristics of the model. In
turn these representations are then easily used to solve the
appropriate control problems formulated below.

In the literature, the dam, storage and production models are
associated with state-dependent queueing systems (e.g.
\cite{Abdel-Hameed 2000}, \cite{Abdel-Hameed Nakhi 1990},
\cite{Bae Kim Lee 2002}, \cite{Bae Kim Lee 2003}, \cite{Boxma
Kaspi Kella Perry 2005}, \cite{Faddy 1974}, \cite{Kaspi Kella
Perry 1996},  \cite{Lam Lou 1987}, \cite{Lee Ahn 1998},
\cite{Zukerman 1977} and others). The model of the present paper
is the following. We assume that units of water, arriving to a
dam, are registered by counter at random instants $t_1$,
$t_2$,\ldots, and interarrival times $\tau_n=t_{n+1}-t_n$ are
mutually independent exponentially distributed random variables
with parameter $\lambda$. Outflow of water is state-dependent as
follows. If the level of water is between $L^{lower}$ and
$L^{upper}$, then an interval between unit departures has the
probability distribution $B_1(x)$. If the level of water increases
above the level $L^{upper}$, then the probability distribution of
an interval between unit departures is $B_2(x)$. It is also
assumed that if the level of water is exactly of the level
$L^{lower}$, then the departure process of water is frozen, and it
is resumed again as soon as the level of water increases the value
$L^{lower}$. It is worth noting that the policies at which the
service rate is changed in dependence of a dam level are of
notable attention in the literature (see \cite{Abdel-Hameed 2000},
\cite{Abdel-Hameed Nakhi 1990}, \cite{Bae Kim Lee 2002}, \cite{Bae
Kim Lee 2003}, \cite{Zukerman 1977} and other papers). However all
of them discuss performance measures associated with an
appropriate upper level of water in a dam
 and, to our best
knowledge, the known results on performance analysis of river
flows structured by lower and upper levels are analytically
difficult and hard for real applications even for simple models
(e.g. see review paper \cite{Phatarfod 1989}). Furthermore, in the
most studies the explicit representations are in terms of the
Laplace-Stieltjes transforms of the initial distributions, and
there is no paper providing asymptotic analysis of large dams.

In terms of queueing theory the problem considered in the paper
can be reformulated as follows. Consider single-server queueing
system where arrival flow of customers is Poisson with rate
$\lambda$, and a service time of a customer depends upon
queue-length as follows. If at the moment of service begun, the
number of customers in the system is not greater than $L$, then
the service time of this customer has the probability distribution
$B_1(x)$. Otherwise, if there are more than $L$ customers in the
system at the moment of service begun, then the probability
distribution function of the service time of this customer is
$B_2(x)$. The analytical results for this queueing system are
known (e.g. Abramov \cite{Abramov 1991}). Notice that the lower
level of dam $L^{lower}$ is equated with an empty queueing system.
Then the dam specification of the problem is characterized by
performance criteria, which in terms of the queueing formulation
looks as follows. Let $q_t$ denote the queue-length in time $t$.
The problem is to choose an output parameter of system minimizing
the functional $J(L)=p_1(L)J_1(L)+p_2(L)J_2(L)$, where
$p_1(L)=\lim_{t\to\infty}\mathbf{P}\{q_t=0\}$,
$p_2(L)=\lim_{t\to\infty}\mathbf{P}\{q_t>L\}$, and $J_1(L)$ and
$J_2(L)$ are the corresponding damage costs proportional to $L$.
 To be precise assume that $J_1(L)=j_1L$ and
$J_2(L)=j_2L$, where $j_1$ and $j_2$ are positive constants.
Assuming that $L\to\infty$ we shall often write $p_1$ and $p_2$
(without argument $L$) rather than $p_1(L)$ and $p_2(L)$. The
argument $L$ will be often omitted in other functions. We shall
feel free to write $J$, $J_1$ and $J_2$ without the argument $L$.

To specify the problem more correctly we assume that the input
parameter $\lambda$, and probability distribution function
$B_2(x)$ are given, while $B_1(x)=B_1(x, C)$ is a family of
probability distributions depending on  parameter $C\geq0$, which
in turn is closely related with the expectation $\int_0^\infty
x\mbox{d}B_1(x)$. Then the output rate, associated with this
probability distribution $B_1(x)$, can be changed so that the
minimum value of the functional is associated with the choice of
this parameter $C$, resulting in the choice of the function
$B_1(x,C)$. The correctness of such formulation and more concrete
clarification of parameter $C$ will be explained in the sequel
(see the formulations of Theorem \ref{thm2}, \ref{thm3},
\ref{thm5}, \ref{thm6}). It is interesting to note, \emph{that the
solution of the above control problem is asymptotically
independent of the explicit form of probability distribution
functions $B_1(x)$ and $B_2(x)$, and only depends on the
expectations}
 $\int_0^\infty
x\mbox{d}B_2(x)<\frac{1}{\lambda}$ \emph{and} $\int_0^\infty
x\mbox{d}B_1(x)$ \emph{as well as} $\int_0^\infty
x^2\mbox{d}B_1(x)$. The details of this dependence will be
explained later.

We use the notation $b_i=\int_0^\infty x\mbox{d}B_i(x)$,
$\rho_i=\lambda b_i$ ($i=1,2$) and assume that $\rho_2<1$. This
assumption is a standard condition of stationarity, ergodicity of
the queue-length process $q_t$ and existence of the limits $p_1$
and $p_2$ (independent of any initial state of the process). In
additional to this assumption we shall also assume the existence
of the third moment: $\rho_{1,k}={\lambda^k}\int_0^\infty
x^k\mbox{d}B_1(x)<\infty$, \ $k=2,3$. The existence of the second
moment is used in Theorem \ref{thm1}. Then the existence of the
third moment for all the specified family of distributions
$B_1(x,C)$ is required in Theorem \ref{thm2}, \ref{thm3},
\ref{thm5}, \ref{thm6}.

The special features of the present paper are as follows.
\smallskip

$\bullet$ We solve the control problem where the performance
criteria takes into account passage of the upper and lower bounds.
The formulation of the problem is not traditional but realistic.
Other similar control problems arising in practice can be flexibly
solved by adapting the method of this paper.

\smallskip
$\bullet$ The presentation of the results are clear and available
for real application.

\smallskip
$\bullet$ The mathematical methods of the paper are not
traditional but clear and easily understandable.

\smallskip
The paper is structured as follows. In Section 2 we discuss the
state-dependent queue-length process and derive representation for
the probabilities $p_1$ and $p_2$. Section 3 contains the results
on asymptotic analysis of probabilities $p_1$ and $p_2$, and the
main result of this asymptotic behaviour is given by Theorem
\ref{thm1}. In Section 4  some additional theorems on asymptotic
behaviour of $p_1$ and $p_2$ are proved, which are then used to
solve the control problem. The main result of this paper, the
solution of control problem is formulated in Section 5. Concluding
remark are given in Section 6.

\section{\bf The state-dependent queue and its characteristics in a busy
period} In this section we discuss the main characteristics of the
state-dependent queueing system described in the introduction. Let
$T_L$, $I_L$ and $\nu_L$ denote correspondingly a busy period, an
idle period and the number of served customers during a busy
period. Let $T_L^{(1)}$, $T_L^{(2)}$ denote the total time during
a busy period when correspondingly $0<q_t\leq L$ and $q_t>L$, and
let $\nu_L^{(1)}$, $\nu_L^{(2)}$ denote correspondingly the total
numbers of served customers during a busy period when
correspondingly $0< q_t\leq L$ and $q_t>L$. We have the following
two obvious equations:
\begin{equation}
\label{2.1}\mathbf{E}T_L=\mathbf{E}T_L^{(1)}+\mathbf{E}T_L^{(2)},
\end{equation}
\begin{equation}
\label{2.2}\mathbf{E}\nu_L=\mathbf{E}\nu_L^{(1)}+\mathbf{E}\nu_L^{(2)}.
\end{equation}
According to the Wald's equation,
\begin{equation}
\label{2.3}\mathbf{E}T_L^{(1)}=b_1\mathbf{E}\nu_L^{(1)},
\end{equation}
and
\begin{equation}
\label{2.4}\mathbf{E}T_L^{(2)}=b_2\mathbf{E}\nu_L^{(2)}.
\end{equation}
Next, the number of arrivals during a busy circle coincides with
the total number of served customers during a busy period. Hence,
applying the Wald's equation again and taking into account
\eqref{2.1}-\eqref{2.4}, we obtain
\begin{equation}\label{2.5}
\begin{aligned}
\lambda\mathbf{E}T_L+\lambda\mathbf{E}I&=\lambda\mathbf{E}T_L+1\\
&=\lambda\mathbf{E}T_L^{(1)}+\lambda\mathbf{E}T_L^{(2)}+1\\
&=\rho_1\mathbf{E}\nu_L^{(1)}+\rho_2\mathbf{E}\nu_L^{(2)}+1\\
&=\mathbf{E}\nu_L^{(1)}+\mathbf{E}\nu_L^{(2)}.
\end{aligned}
\end{equation}
From \eqref{2.5} we have the equation
\begin{equation}
\label{2.6}\mathbf{E}\nu_L^{(2)}=\frac{1}{1-\rho_2}-\frac{1-\rho_1}{1-\rho_2}\mathbf{E}\nu_L^{(1)},
\end{equation}
expressing $\mathbf{E}\nu_L^{(2)}$ via $\mathbf{E}\nu_L^{(1)}$.
For example, if $\rho_1=1$, then
$\mathbf{E}\nu_L^{(2)}=\frac{1}{1-\rho_2}$ for any $L$.

The similar equation holds also for $\mathbf{E}T_L^{(2)}$. Namely,
from \eqref{2.4} and \eqref{2.6} we obtain
\begin{equation}
\label{2.7}\mathbf{E}T_L^{(2)}=\frac{\rho_2}{\lambda(1-\rho_2)}
-\frac{\rho_2(1-\rho_1)}{\lambda(1-\rho_2)}\mathbf{E}T_L^{(1)}.
\end{equation}
Equations \eqref{2.6} and \eqref{2.7} enables us to obtain the
stationary probabilities $p_1$ and $p_2$. Applying the renewal
reward theorem (e.g. Ross \cite{Ross 1983}, p. 78) and
consequently \eqref{2.5} and \eqref{2.6}, for $p_1$ we obtain:
\begin{equation}\label{2.8}
\begin{aligned}
p_1&=\frac{\mathbf{E}I}{\mathbf{E}T_L^{(1)}+\mathbf{E}T_L^{(2)}+\mathbf{E}I}\\
&=\frac{1}{\mathbf{E}\nu_L^{(1)}+\mathbf{E}\nu_L^{(2)}}\\
&=\frac{1-\rho_2}{1+(\rho_1-\rho_2)\mathbf{E}\nu_L^{(1)}}.\\
\end{aligned}
\end{equation}
Analogously,
\begin{equation}
\label{2.9}
\begin{aligned}
p_2&=\frac{\mathbf{E}T_L^{(2)}}{\mathbf{E}T_L^{(1)}+\mathbf{E}T_L^{(2)}+\mathbf{E}I}\\
&=\frac{\rho_2\mathbf{E}\nu_L^{(2)}}{\mathbf{E}\nu_L^{(1)}+\mathbf{E}\nu_L^{(2)}}\\
&=\frac{\rho_2+\rho_2(\rho_1-1)\mathbf{E}\nu_L^{(1)}}{1+(\rho_1-\rho_2)\mathbf{E}\nu_L^{(1)}}.\\
\end{aligned}
\end{equation}

\section{\bf Asymptotic analysis of $p_1$ and $p_2$ as $L$ increases
to infinity} By sample path analysis and the property of the lack
of memory of exponential distribution it follows that the random
variable $\nu_L^{(1)}$ coincides in distribution with the number
of served customers during a busy period of the $M/GI/1/L$
queueing system (the parameter $L$ denotes the number of customers
in the system excluding the customer in the server). Specifically,
we use the fact that during a busy period the number of times
 service begun when the number of customers in the
system does not exceed $L$, coincides with the number of arrivals
when the number of customers in the system does not exceed $L+1$.
We also use the fact that the residual interarrival time after a
service completion has exponential distribution with parameter
$\lambda$.

Therefore the known results of the $M/GI/1/L$ queueing system can
be used.

It is known (e.g. \cite{Abramov 1991}, \cite{Abramov 1997}) that
$\mathbf{E}\nu_L^{(1)}$ is determined by the convolution type
recurrence relation
\begin{equation*}
\mathbf{E}\nu_L^{(1)}=\sum_{j=0}^{L}\mathbf{E}\nu_{L-j+1}^{(1)}\int_0^\infty
\mbox{e}^{-\lambda x}\frac{(\lambda x)^j}{j!}\mbox{d}B_1(x), \ \ \
\mathbf{E}\nu_0^{(1)}=1,
\end{equation*}
where $\mathbf{E}\nu_{n}^{(1)}$ denotes the expectation of the
number of served customers during a busy period of $M/GI/1/n$
queue ($n=1,2,\ldots$).

The probabilities $p_1(L)$ and $p_2(L)$ are expressed explicitly
via $\mathbf{E}\nu_L^{(1)}$, and their asymptotic behavior as
$L\to\infty$ can be obtained from the following known results.

Let $Q_0> 0$ be an arbitrary real number, and for $n\geq 0$
\begin{equation*}
Q_n=\sum_{j=0}^nr_jQ_{n-j+1},
\end{equation*}
where $r_0>0$, $r_j\geq0$, and $r_0+r_1+\ldots=1$. Let
$r(z)=\sum_{j=0}^\infty r_jz^j$, $|z|\leq1$ be a generating
function, and let $\gamma_m=\lim_{z\uparrow 1}r^{(m)}(z)$, where
$r^{(m)}(z)$ is the $m$th derivative of $r(z)$.

Notice, that the sequence $\{Q_n\}$ is an increasing sequence, and
\begin{equation}
\label{3.0}\sum_{n=0}^\infty Q_nz^n=\frac{Q_0r(z)}{z-r(z)}
\end{equation}
(see \cite{Postnikov 1979}, \cite{Postnikov 1980},  Sect. 25 and
\cite{Takacs 1967}).

\begin{lem}
\label{lem1} (Tak\'acs \cite{Takacs 1967}, p. 22-23). If
$\gamma_1<1$, then
\begin{equation*}
\lim_{n\to\infty}Q_n=\frac{Q_0}{1-\gamma_1}.
\end{equation*}
If $\gamma_1=1$ and $\gamma_2<\infty$, then
\begin{equation*}
\lim_{n\to\infty}\frac{Q_n}{n}=\frac{2Q_0}{\gamma_2}.
\end{equation*}
If $\gamma_1>1$, then
\begin{equation*}
\lim_{n\to\infty}\left[Q_n-\frac{Q_0}{\sigma^n(1-r^\prime(\sigma))}\right]
=\frac{1}{1-\gamma_1},
\end{equation*}
where $\sigma$ is the least in absolute value root of functional
equation $z=r(z)$.
\end{lem}



From this lemma we have the asymptotic results for the
probabilities $p_1$ and $p_2$.

For $\Re(s)\geq0$ denote by $\widehat B_1(s)$ the
Laplace-Stieltjes transform of $B_1(x)$. We have the following
theorem.

\begin{thm}
\label{thm1} If $\rho_1<1$, then
\begin{eqnarray}
\label{3.1}\lim_{L\to\infty}p_1(L)&=&1-\rho_1,\\
\label{3.2}\lim_{L\to\infty}p_2(L)&=&0.
\end{eqnarray}
If $\rho_1=1$, then
\begin{eqnarray}
\label{3.3}\lim_{L\to\infty}Lp_1(L)&=&\frac{\rho_{1,2}}{2},\\
\label{3.4}\lim_{L\to\infty}Lp_2(L)&=&\frac{\rho_2}{1-\rho_2}\cdot\frac{\rho_{1,2}}{2}.
\end{eqnarray}
If $\rho_1>1$, then
\begin{equation}\label{3.5}
\lim_{L\to\infty}\frac{p_1(L)}{\varphi^L}=\frac{(1-\rho_2)[1+\lambda
\widehat B_1^\prime(\lambda-\lambda\varphi)]}{\rho_1-\rho_2},
\end{equation}
where $\varphi$ is the least in absolute value root of the
functional equation $z=\widehat B_1(\lambda-\lambda z)$, and
\begin{equation}\label{3.6}
\lim_{L\to\infty}p_2(L)=\frac{\rho_2(\rho_1-1)}{\rho_1-\rho_2}.
\end{equation}
\end{thm}

\begin{proof}
The proof of this theorem follows by application of Lemma
\ref{lem1}. Straightforward application of the aforementioned
lemma to the recurrence relation for $\mathbf{E}\nu_L^{(1)}$
yields the following.

If $\rho_1<1$, then
\begin{equation}\label{3.7}
\lim_{L\to\infty}\mathbf{E}\nu_L^{(1)}=\frac{1}{1-\rho_1}.
\end{equation}
If $\rho_1=1$, then
\begin{equation}\label{3.8}
\lim_{L\to\infty}\frac{\mathbf{E}\nu_L^{(1)}}{L}=\frac{2}{\rho_{1,2}}.
\end{equation}
If $\rho_1>1$, then
\begin{equation}\label{3.9}
\lim_{L\to\infty}\left[\mathbf{E}\nu_L^{(1)}-\frac{1}{\varphi^{L}(1+\lambda\widehat
B_1^\prime(\lambda-\lambda z))}\right]=\frac{1}{1-\rho_1}.
\end{equation}
Substituting \eqref{3.7}-\eqref{3.9} for the limits in \eqref{2.8}
and \eqref{2.9} correspondingly to the cases in the formulation of
the theorem finishes the proof. \end{proof}

\section{\bf Further asymptotic analysis of $p_1$ and $p_2$}
Let us first discuss the statements of Theorem \ref{thm1}. Under
the assumption $\rho_1<1$ we have \eqref{3.1} and \eqref{3.2}. The
probability $p_1$ is positive in limit while the probability $p_2$
vanishes. Under the assumption $\rho_1>1$ we have \eqref{3.5} and
\eqref{3.6}. According to these relations the probability $p_1$
vanishes while $p_2$ is positive in limit. This means that if both
$J_1$ and $J_2$ are large positive values proportional to $L$,
then the functional $J$ will take the value proportional to a
large parameter $L$ too. Specifically, in the case $\rho_1<1$ for
this value we have $J\approx(1-\rho_1)J_1$, and in the case
$\rho_1>1$ we have
$J\approx\frac{\rho_2(\rho_1-1)}{\rho_1-\rho_2}J_2$.

In the case $\rho_1=1$ both $p_1$ and $p_2$ vanish with the rate
$L^{-1}$, and therefore $J$ converges to the limit as
$L\to\infty$. Thus, the case $\rho_1=1$
 is a possible solution of the control problem, while the cases $\rho_1<1$
and $\rho_1>1$ are irrelevant. Specifically, for $J=J(L)$ we
obtain the following:
\begin{equation}
\label{4.1}
\lim_{L\to\infty}J(L)=j_1\frac{\rho_{1,2}}{2}+j_2\frac{\rho_2}{1-\rho_2}\cdot\frac{\rho_{1,2}}{2}.
\end{equation}

In order to find now the optimal solution consider the following
two cases: (i) \ $\rho_1=1+\delta$ and (ii) \ $\rho_1=1-\delta$,
where in these both  cases $\delta\to 0$ as $L\to\infty$.
\smallskip


In case (i) we have the following two theorems.

\begin{thm}
\label{thm2} Assume that $\rho_1=1+\delta$, $\delta>0$, and
$L\delta\to C>0$, as $\delta\to 0$ and $L\to\infty$. Assume that
$\rho_{1,3}=\rho_{1,3}(\delta)$ is a bounded function of parameter
$\delta$, for all $0\leq\delta<1$ and there exists
$\widetilde\rho_{1,2}=\lim_{\delta\to 0}\rho_{1,2}(\delta)$. Then,
\begin{eqnarray}
\label{4.2}p_1&=&\frac{\delta}{\mbox{e}^{2C/\widetilde\rho_{1,2}}-1}+o(\delta),\\
\label{4.3}p_2&=&\frac{\delta\rho_2\mbox{e}^{2C/\widetilde\rho_{1,2}}}
{(1-\rho_2)(\mbox{e}^{2C/\widetilde\rho_{1,2}}-1)}+o(\delta).
\end{eqnarray}

\end{thm}

\begin{proof}
The proof of this theorem is similar to that of Theorem 3.4 of
\cite{Abramov 2002} and Theorem 4.4 of \cite{Abramov 2004}. Under
the conditions of the theorem the following expansion was shown in
Subhankulov \cite{Subhankulov 1976}, p. 326:
\begin{equation}
\label{4.4}\varphi=1-\frac{2\delta}{\widetilde\rho_{1,2}}+O(\delta^2).
\end{equation}
Then, by virtue of \eqref{4.4} after some algebra we have:
\begin{equation}
\label{4.5}1+\lambda\widehat
B^\prime(\lambda-\lambda\varphi)=\delta+O(\delta^2).
\end{equation}
Substituting \eqref{4.4} and \eqref{4.5} for \eqref{3.9} we
obtain:
\begin{equation}
\label{4.6}\mathbf{E}\nu_L^{(1)}=\frac{\mbox{e}^{2C/\widetilde\rho_{1,2}}-1}{\delta}+O(1).
\end{equation}
From \eqref{4.6} and \eqref{2.8} and \eqref{2.9} we finally obtain
the statement of the theorem.
\end{proof}

\begin{thm}
\label{thm3}Under the conditions of Theorem \ref{thm2} assume that
$C=0$. Then,
\begin{eqnarray}
\label{4.7}\lim_{L\to\infty}Lp_1(L)&=&\frac{\rho_{1,2}}{2},\\
\label{4.8}\lim_{L\to\infty}Lp_2(L)&=&\frac{\rho_2}{1-\rho_2}\cdot\frac{\rho_{1,2}}{2}.
\end{eqnarray}
\end{thm}
\begin{proof}
The statement of the theorem follows by expanding the main terms
of asymptotic relations of \eqref{4.2} and \eqref{4.3} for small
$C$.
\end{proof}

Notice, that \eqref{4.7} and \eqref{4.8} coincide with \eqref{3.3}
and \eqref{3.4} correspondingly.

\smallskip
In case (ii) we have the following.

\begin{thm}
\label{thm5}Assume that $\rho_1=1-\delta$, $\delta>0$, and
$L\delta\to C>0$, as $\delta\to 0$ and $L\to\infty$. Assume that
$\rho_{1,3}=\rho_{1,3}(\delta)$ is a bounded function of parameter
$\delta$, for all $0\leq\delta<1$ and there exists
$\widetilde\rho_{1,2}=\lim_{\delta\to 0}\rho_{1,2}(\delta)$. Then,
\begin{eqnarray}
\label{4.9}p_1&=&\delta\mbox{e}^{\widetilde\rho_{1,2}/2C}+o(\delta),\\
\label{4.10}p_2&=&\delta\cdot\frac{\rho_2}{1-\rho_2}\left(\mbox{e}^{\widetilde\rho_{1,2}
/2C}-1\right)+o(\delta).
\end{eqnarray}

\end{thm}

\begin{proof}
From \eqref{3.0} we have
\begin{equation*}\label{4.11}
\sum_{n=0}^\infty \mathbf{E}\nu_n^{(1)}z^{n}=\frac{\widehat
B_1(\lambda-\lambda z)}{\widehat B_1(\lambda-\lambda z)-z}.
\end{equation*}
The sequence $\{\mathbf{E}\nu_n^{(1)}\}$ is an increasing
sequence, and in the case $\rho_1=1$ from the Tauberian theorem of
Hardy-Littlewood (e.g. \cite{Postnikov 1979}, \cite{Postnikov
1980}, \cite{Subhankulov 1976}, \cite{Sznajder Filar 1992},
\cite{Takacs 1967}) we obtain:
\begin{equation*}
\lim_{L\to\infty}\frac{\mathbf{E}\nu_L^{(1)}}{L}=\lim_{z\uparrow
1}(1-z)^2\frac{\widehat B_1(\lambda-\lambda z)}{\widehat
B_1(\lambda-\lambda z)-z}.
\end{equation*}
(It is not difficult to check that then \eqref{3.8} follows.) Then
in the case where $\rho_1=1-\delta$, and $L\delta\to C$ as
$L\to\infty$, according to the same Tauberian theorem of Hardy and
Littlewood, asymptotic behaviour of $\mathbf{E}\nu_L^{(1)}$ can be
found from the asymptotic expansion
\begin{equation}\label{4.12}
(1-z)\cdot\frac{\widehat B_1\left(\lambda-\lambda
z\right)}{\widehat B_1\left(\lambda-\lambda z\right)-z},
\end{equation}
as $z\uparrow 1$.

By the Taylor expansion of the denominator of \eqref{4.12} we
obtain:
\begin{equation}\label{4.13}
\begin{aligned}
\frac{1-z}{\widehat B_1\left(\lambda-\lambda z\right)-z}&\asymp
\frac{1-z}{1-z-\rho_1(1-z)+\dfrac{\widetilde\rho_{1,2}}{2}(1-z)^2+O\left((1-z)^3\right)}\\
&\asymp\frac{1}{\delta+\dfrac{\widetilde\rho_{1,2}}{2}(1-z)+O\left[(1-z)^2\right]}\\
&\asymp\frac{1}{\delta\left(1+\dfrac{\widetilde\rho_{1,2}}{2\delta}(1-z)+O\left[(1-z)^2
\right]\right)}\\
&\asymp\frac{1}{\delta\exp\left(\dfrac{\widetilde\rho_{1,2}}{2\delta}(1-z)\right)}\cdot[1+o(1)].
\end{aligned}
\end{equation}
Therefore, assuming that $z=\frac{L-1}{L}\to 1$ as $L\to\infty$,
from \eqref{4.13} we obtain the asymptotic behaviour of
$\mathbf{E}\nu_L^{(1)}$ as $L\to\infty$. We have:
\begin{equation}
\label{4.14}\mathbf{E}\nu_L^{(1)}=\frac{1}{\delta\mbox{e}^{\widetilde\rho_{1,2}/2C}}
\cdot [1+o(1)].
\end{equation}
Now, substituting \eqref{4.14} for \eqref{2.8} and \eqref{2.9} we
obtain the desired statements of the theorem.
\end{proof}

\begin{thm}
\label{thm6}Under the conditions of Theorem \ref{thm5} assume that
$C=0$. Then we obtain \eqref{4.7} and \eqref{4.8}.
\end{thm}

\begin{proof}
The statement of the theorem follows by expanding the main terms
of asymptotic relations of \eqref{4.9} and \eqref{4.10} for small
$C$.
\end{proof}

\section{\bf Solution of the control problem}
In this section we formulate the theorem characterizing the
solution of control problem.

For $J=J(L)$ we have the following limiting relation
\begin{equation}
\begin{aligned}
\label{5.1}\lim_{L\to\infty}J(L)&=\lim_{L\to\infty}\left[p_1(L)J_1(L)+p_2(L)J_2(L)\right]\\
&=j_1\lim_{L\to\infty}Lp_1(L)+j_2\lim_{L\to\infty}Lp_2(L).
\end{aligned}
\end{equation}
Substituting \eqref{4.1} and \eqref{4.2} for the right-hand side
of \eqref{5.1} and taking into account that $L\delta\to C$, we
obtain:
\begin{equation}
\begin{aligned}
\label{5.2}J^{upper}&=\lim_{L\to\infty}J(L)\\
&=\left[j_1\frac{1}{\mbox{e}^{2C/\widetilde\rho_{1,2}}-1}+
j_2\frac{\rho_2\mbox{e}^{2C/\widetilde\rho_{1,2}}}
{(1-\rho_2)(\mbox{e}^{2C/\widetilde\rho_{1,2}}-1)}\right]\lim_{L\to\infty}L\delta\\
&=C\left[j_1\frac{1}{\mbox{e}^{2C/\widetilde\rho_{1,2}}-1}+
j_2\frac{\rho_2\mbox{e}^{2C/\widetilde\rho_{1,2}}}
{(1-\rho_2)(\mbox{e}^{2C/\widetilde\rho_{1,2}}-1)}\right].
\end{aligned}
\end{equation}
Substituting \eqref{4.9} and \eqref{4.10} for the right-hand side
of \eqref{5.1} and taking into account that $L\delta\to C$, we in
turn obtain:
\begin{equation}
\begin{aligned}
\label{5.3}J^{lower}
&=C\left[j_1\mbox{e}^{\widetilde\rho_{1,2}/2C}+
j_2\frac{\rho_2}{1-\rho_2}\left(\mbox{e}^{\widetilde\rho_{1,2}
/2C}-1\right)\right].
\end{aligned}
\end{equation}
Let us now study the functionals $J^{upper}$ and $J^{lower}$ given
by \eqref{5.2} and \eqref{5.3}. Observing \eqref{5.2}, notice that
there contain the constants $j_1$, $j_2$ and $\rho_2$ in
\eqref{5.2}. Let us assume that these constants are given such
that
\begin{equation}
\label{5.4}j_1=j_2\cdot\frac{\rho_2}{1-\rho_2}.
\end{equation}
Then $C=0$ is the point of $\min$ of the functional $J^{upper}$.
Indeed, in this case
\begin{equation}
\label{5.5}
\begin{aligned}
J^{upper}&=j_1C\left[\frac{1}{\mbox{e}^{2C/\widetilde\rho_{1,2}}-1}+
\frac{\mbox{e}^{2C/\widetilde\rho_{1,2}}}{\mbox{e}^{2C/\widetilde\rho_{1,2}}-1}\right]\\
&=j_1C\left[\frac{1}{\mbox{e}^{2C/\widetilde\rho_{1,2}}-1}+
\frac{\left(\mbox{e}^{2C/\widetilde\rho_{1,2}}-1\right)+1}
{\mbox{e}^{2C/\widetilde\rho_{1,2}}-1}\right]\\
&=j_1C\left[1+\frac{2}{\mbox{e}^{2C/\widetilde\rho_{1,2}}-1}\right].\\
\end{aligned}
\end{equation}
Therefore in point $C=0$ we have $\lim_{C\to
0}J^{upper}=j_1\widetilde\rho_{1,2}$, and in the right side of the
point $C=0$ the function $J^{upper}$ is increasing in $C$. Hence
\eqref{5.4} is the condition for $C=0$.

Next,
\begin{equation}\label{5.6}
\left[\frac{C}{\mbox{e}^{2C/\widetilde\rho_{1,2}}-1}\right]_C^\prime=
\frac{\mbox{e}^{2C/\widetilde\rho_{1,2}}-1-\dfrac{2C^2}{\widetilde\rho_{1,2}}
\cdot\mbox{e}^{2C/\widetilde\rho_{1,2}}}
{(\mbox{e}^{2C/\widetilde\rho_{1,2}}-1)^2},
\end{equation}
and
\begin{equation}\label{5.7}
\begin{aligned}
\left[\frac{C\mbox{e}^{2C/\widetilde\rho_{1,2}}}{\mbox{e}^{2C/\widetilde\rho_{1,2}}-1}
\right]_C^\prime&=\left[\frac{C}{\mbox{e}^{2C/\widetilde\rho_{1,2}}-1}\right]_C^\prime
\mbox{e}^{2C/\widetilde\rho_{1,2}}\\
&+\frac{2C}{\widetilde\rho_{1,2}}\left[\frac{C\mbox{e}^{2C/\widetilde\rho_{1,2}}}{\mbox{e}^{2C/\widetilde\rho_{1,2}}-1}
\right]\mbox{e}^{2C/\widetilde\rho_{1,2}}.
\end{aligned}
\end{equation}

Clearly that \eqref{5.7} is not smaller that \eqref{5.6}, and they
are equal when $C=0$.

Therefore, if the right-hand side of \eqref{5.4} is greater than
that left-hand side of \eqref{5.4}, then $C=0$ remains to be the
value minimizing the functional $J^{upper}$.
The similar result holds for functional $J^{lower}$ given in
\eqref{5.3}. Specifically, if the right-hand side of \eqref{5.4}
is not greater than the left-hand side of \eqref{5.4}, then $C=0$
remains to be the value minimizing the functional $J^{lower}$.

Thus, the solution of control problem is given by the following
theorem.

\begin{thm}
\label{thm4}If the parameters $\lambda$ and $\rho_2$ are given,
then the optimal solution of the control problem is the following.

$\bullet$ If
\begin{equation*}
j_1=\frac{\rho_2}{1-\rho_2}j_2,
\end{equation*}
then the optimal solution of the control problem is achieved for
$\rho_1=1$.

$\bullet$ If
\begin{equation*}
j_1>\frac{\rho_2}{1-\rho_2}j_2,
\end{equation*}
then the optimal solution of the control problem is a minimization
of the functional $J^{upper}$. The optimal solution is achieved
for $\rho_1=1+\delta$, $\delta(L)$ is a small positive parameter,
and $L\delta(L)\to C$. $C$ is the nonnegative parameter minimizing
\eqref{5.2}.

$\bullet$ If
\begin{equation*}
j_1<\frac{\rho_2}{1-\rho_2}j_2,
\end{equation*}
then the optimal solution of the control problem is a minimization
of the functional $J^{lower}$. The optimal solution is achieved
for $\rho_1=1-\delta$, $\delta(L)$ is a small positive parameter,
and $L\delta(L)\to C$. $C$ is the nonnegative parameter minimizing
\eqref{5.3}.
\end{thm}

\section{\bf Concluding remarks}
In this paper we posed and solved a control problem for a large
dam. The main specification of the problem is that the performance
criteria takes into account passage the lower and upper bounds.
The solution of the control problem is asymptotically independent
of the explicit form of probability distribution functions
$B_1(x)$ and $B_2(x)$, and under the assumption that the
parameters $\lambda$ and $\rho_2$ are given, in dependence of a
performance criteria the parameter $\rho_1$ must have one of the
forms: $\rho_1=1$, $\rho_1=1+\delta(L)$, or $\rho_1=1-\delta(L)$
where $\delta(L)>0$, and as $L\to\infty$, $\delta(L)$ vanishes and
$L\delta(L)\to C$.

\section*{\bf Acknowledgement}
The research was supported by Australian Research Council grant
No. DP0771338.

\bibliographystyle{amsplain}

\begin{thebibliography}{10}
\bibitem{Abdel-Hameed 2000} \textsc{Abdel-Hameed, M.S.} (2000).
Optimal control of a dam using $P_{\lambda,\tau}^M$ policies and
penalty cost when the input process is a compound Poisson process
with positive drift. \emph{Journal of Applied Probability}, 37,
406-416.

\bibitem{Abdel-Hameed Nakhi 1990}\textsc{Abdel-Hameed, M.S. and
Nakhi, Y.} (1990). Optimal control of a finite dam using
$P_{\lambda,\tau}^M$ policies and penalty cost: total discounted
and long-run average cases. \emph{Journal of Applied Probability},
27, 888-898.

\bibitem{Abramov 1991} \textsc{Abramov, V.M.} (1991).
\emph{Investigation of a Queueing System with Service Depending on
a Queue-Length}. Donish, Dushanbe, Tadzhikistan. (Russian.)

\bibitem{Abramov 1997} \textsc{Abramov, V.M.} (1997). On a
property of a refusals stream. \emph{Journal of Applied
Probability}, 37, 800-805.

\bibitem{Abramov 2002} \textsc{Abramov, V.M.} (2002). Asymptotic
analysis of the $GI/M/1/n$ queueing system as $n$ increases to
infinity. \emph{Annals of Operations Research}, 112, 35-41.

\bibitem{Abramov 2004} \textsc{Abramov, V.M.} (2004). Asymptotic
behavior of the number of lost messages. \emph{SIAM Journal on
Applied Mathematics} 64 (3) 746-761.

\bibitem{Bae Kim Lee 2002}\textsc{Bae, J., Kim, S. and Lee, E.Y.}
(2002). A $P_\lambda^M$ policy for an M/G/1 queueing system.
\emph{Applied Mathematical Modelling}, 26, 929-939.

\bibitem{Bae Kim Lee 2003}\textsc{Bae, J., Kim, S. and Lee, E.Y.}
(2003). Average cost under the $P_{\lambda,\tau}^M$ policy in a
finite dam with compound Poisson inputs. \emph{Journal of Applied
Probability}, 40, 519-526.

\bibitem{Boxma Kaspi Kella Perry 2005}\textsc{Boxma, O., Kaspi,
H., Kella, O. and Perry, D.} (2005). On/off storage systems with
state-dependent input, output, and switching rates.
\emph{Probability in the Engineering and Informational Sciences},
19, 1-14.

\bibitem{Faddy 1974}\textsc{Faddy, M.J.} (1974).
Optimal control of finite dams: discrete (2-stage) output
procedure. \emph{Journal of Applied Probability}, 11, 111-121.

\bibitem{Kaspi Kella Perry 1996}\textsc{Kaspi, H., Kella, O.,
Perry, D.} (1996). Dam processes with state-dependent batch sizes,
and intermittent production processes with state-dependent rates.
\emph{Queueing Systems}, 24, 37-57.

\bibitem{Lam Lou 1987}\textsc{Lam, Y. and Lou, J.H.} (1987).
Optimal control for a finite dam. \emph{Journal of Applied
Probability}, 24, 196-199.

\bibitem{Lee Ahn 1998}\textsc{Lee, E.Y. and Ahn, S.K.} (1998).
$P_\tau^M$ policy for a dam with input formed by a compound
Poisson process. \emph{Journal of Applied Probability}, 35,
482-488.

\bibitem{Phatarfod 1989}\textsc{Phatarfod, R.M.} (1989). Riverflow
and reservoir storage models. \emph{Mathematical and Computer
Modelling}, 12, 1057-1077.

\bibitem{Postnikov 1979}\textsc{Postnikov, A.G.} (1979).
\emph{Tauberian Theory and its Application}. Trudy Mat. Inst.
Steklov, (2) 144. (Russian).

\bibitem{Postnikov 1980}\textsc{Postnikov, A.G.} (1980).
\emph{Tauberian Theory and its Application}. Proc. Steklov Math.
Inst., (2) 144. (AMS transl. from Russian.)


\bibitem{Ross 1983} \textsc{Ross, S.M.} (1983). \emph{Stochastic
Processes}, John Wiley, New York.

\bibitem{Subhankulov 1976}\textsc{Subhankulov, M.A.} (1976).
\emph{Tauberian Theorems with Remainder}. Nauka, Moscow.
(Russian.)

\bibitem{Sznajder Filar 1992}\textsc{Sznajder, R. and Filar, J.A.}
(1992). Some comments on a theorem of Hardy and Littlewood.
\emph{Journal of Optimization Theory and Applications}, 75,
201-208.

\bibitem{Takacs 1967}\textsc{Tak\'acs, L.} (1967).
\emph{Combinatorial Methods in the Theory of Stochastic
Processes}, John Wiley, New York.

\bibitem{Zukerman 1977}\textsc{Zukerman, D.} (1977). Two-stage
output procedure of a finite dam. \emph{Journal of Applied
Probability}, 14, 421-425.


\end{thebibliography}

\end{document}